\documentclass{article}

\usepackage{latexsym}
\usepackage[T1]{fontenc}
\usepackage{amssymb,amsmath}
\usepackage{amsthm}
\usepackage{amsrefs}
\usepackage{enumerate}
\usepackage[all]{xy}
\usepackage[OT2,T1]{fontenc}
\usepackage{graphicx}
\DeclareSymbolFont{cyrletters}{OT2}{wncyr}{m}{n}
\DeclareMathSymbol{\Shaa}{\mathalpha}{cyrletters}{"58}

\newcommand{\cls}{\mathcal{H}_{c}}
\newcommand{\opn}{\mathcal{H}_{o}}

\newcommand{\ocha}{\Lambda^c \cls \otimes T^c \opn}

\newcommand{\geqs}{\geqslant}
\newcommand{\leqs}{\leqslant}

\newcommand{\mquad}{\!\!\!\!\!\!}

\newcommand{\rar}{\rightarrow}

\newtheorem{thm}{Theorem}
\newtheorem*{thm*}{Theorem}
\newtheorem{defi}{Definition}
\newtheorem{obs}{Observation}

\newtheorem{prop}{Proposition}
\newtheorem*{prop*}{Proposition}

\newtheorem*{prob*}{Problem}

\newtheorem*{obs*}{Observation}

\begin{document}
\title{On the Coalgebra Description of OCHA}
\author{Eduardo Hoefel\footnote{Supported by CNPq-Brazil grant
        140353/02 (Ph.D. studies at Unicamp) 
        and grant 201064/04 (visiting student at UPenn).} 
\\ {\small hoefel@ufpr.br} }
\maketitle

\begin{abstract}
  \noindent OCHA is the homotopy algebra of open-closed strings. 
   It can be defined as a sequence of multilinear operations on a pair of DG spaces 
   satisfying certain relations which include the $L_\infty$ relations in one 
   space and the $A_\infty$ relations in the other. 
   In this paper we show that the OCHA structure
   is intrinsic to the tensor product of the symmetric and tensor coalgebras.
   We also show how an OCHA can be obtained from $A_\infty$-extesions and define 
   the {\it universal enveloping} $A_\infty$-algebra of an OCHA as an $A_\infty$-extension
   of the universal enveloping of its $L_\infty$ part by its $A_\infty$ part. 
\end{abstract}

\tableofcontents

\section{Introduction} 

Inspired by Zwiebach's classical open-closed string field theory, Kajiura and Stasheff
introduced {\bf O}pen-{\bf C}losed {\bf H}omotopy {\bf A}lgebras \cite{KS06a}. OCHAs were 
presented in three equivalent ways: as a sequence of {\it multilinear operations} satisfying certain conditions; 
as an {\it algebra over a DG Operad} and as a {\it coderivation differential} on a certain coalgebra.

Let $(\cls, \opn)$ be a pair of DG spaces. 
According to the {\it multilinear operations} description, an OCHA structure on 
the pair $(\cls, \opn)$ consists of two sequences of multilinear 
maps: $l_n : \cls^{\land n} \rar \cls$, $n \geqs 1$ and 
$n_{p,q} : \cls^{\land p} \otimes \opn^{\otimes q} \rar \opn$, $p+q \geqs 1$, satisfying certain 
compatbility conditions. 
%
%
%After lifting the maps as coderivations (see 
%formulas (\ref{lift1}) and (\ref{lift2}) on the appendix), the compatbility condition 
%is equivalent to the condition $(\mathfrak{l} + \mathfrak{n})^2 = 0$, where 
%$\mathfrak{l}$ and $\mathfrak{n}$ are the coderivations we obtain after the lifting.
%
The compatibility conditions say that the maps $\{ l_n \}$ define an 
$L_\infty$ structure on $\mathcal{H}_c$, the maps $\{ n_{0,q} \}$ define an $A_\infty$ structure 
on $\mathcal{H}_o$ and the maps $\{ n_{p,q} \}$ for $p,q \geqs 1$ provide the strucutre 
of an $A_\infty$-algebra over an $L_\infty$-algebra on $\mathcal{H}_o$ (i.e., the strongly homotopy version 
of the action by derivations of a Lie algebra on an associative algebra, see \cite{KS06a}). It remais to understand
mathematically the terms of the form $n_{p,0} : \cls^{\land p} \rar \opn$. 
Providing one possible {\it mathematical} interpratation for those maps 
is one of the concerns of the present paper. On the other hand, their physical meaning is 
well recognized.
According to \cite{KS06a}, the maps $n_{p,0}$ originated from 
the operations of {\it opening closed strings into open ones}. 

The {\it operadic description} consists of providing a differential graded operad whose algebras (or representations)
are precisely those pairs of DG spaces endowed with the structure of an OCHA. 
Kajiura and Stasheff defined that operad using the language of trees in \cite{KS06b} and discussed
the geometry behind that operad in \cite{KS06c}. 
That geometrical description was further studied in \cite{Hoe06a} and used 
to give the OCHA operad a description in terms of minimal resolutions involving the Swiss-Cheese operad. 
There we prove that the OCHA operad is the operad defined by 
the first row of the $E^1$ term of the spectral sequence of the compactified configuration space of points on the 
closed disc. 

As for the {\it coderivation description} (or {\it coalgebra description} of OCHA), let us consider
the multlinear maps: $l_n : \cls^{\land n} \rar \cls$ and 
$n_{p,q} : \cls^{\land p} \otimes \opn^{\otimes q} \rar \opn$ satisfying the above mentioned
compatbility conditions. Kajiura and Stasheff showed that, after lifting those maps as coderivations 
$\mathfrak{l}$ and $\mathfrak{n}$ on the coalgebra $\Lambda^c \cls \otimes T^c \opn$, 
the compatbility condition is equivalent to the condition $(\mathfrak{l} + \mathfrak{n})^2 = 0$.
That coalgebra description is the subject of the present paper. We will show that {\it any} coderivation 
on $\Lambda^c \cls \otimes T^c \opn$ has the form $\mathfrak{l} + \mathfrak{n}$, thus showing
that the OCHA structure is {\it intrinsic} to the tensor product of the symmetric and tensor coalgebras.

\subsection*{Main Results}
We now describe the main results of this paper. The notation used will be explained at the end of this introduction.
Let us assume a fixed field $k$ of characteristic zero.
Consider a vector space $E$ with a splitting $E = A \oplus B$ in the category of vector spaces.
Let $M = \{ m_k : E^{\otimes k} \to E \}$ be a family of multilinear maps on $E$. We say that $M$ 
satisfies the OCHA constraint with respect to the spliting $E = A \oplus B$ if it has only components 
of the form $M^A_{p,q} : A^{\otimes p} \otimes B^{\otimes p} \to A$ and 
$M^B_n : B^{\otimes n} \to B$. We say that a coderivation satisfies the constraint if it is 
obtained by lifting a family of maps that satisfy the constraint. More details are given  
in Section \ref{coderivations}, where we prove the following result for $\cls \oplus \opn$.
\begin{prop*}
Any coderivation on the coalgebra $\Lambda^c \cls \otimes T^c \opn$ satisfies the OCHA constraint.
%with respect to the spliting $\cls \oplus \opn$.
\end{prop*}
This fact means that the OCHA constraint, whose geometrical and physical meaning is discussed in 
\cite{KS06c,Hoe06a}, is in fact {\it intrinsic} to the coalgebra $\Lambda^c \cls \otimes T^c \opn$. 
As a consequence, we have the following theorem. 
\begin{thm*}
An OCHA structure on the pair of spaces $(\cls,\opn)$ is equivalent to a degree one coderivation differential 
$\mathcal{D} \in {\rm Coder}(\Lambda^c \cls \otimes T^c \opn)$, $\mathcal{D}^2 = 0$.
\end{thm*}

Section \ref{coshuff} involves a coalgebra map 
$$\Xi : \Lambda^c \cls \otimes T^c \opn \to T^c(\cls \oplus \opn)$$
given by symmetrization and shuffling.
An $A_\infty$-extension is a short exact sequence of $A_\infty$-algebras $0 \to A \to E \to B \to 0$, 
where each map is a linear $A_\infty$-morphism.
The following theorem is proven in Section 4:
\begin{thm*}
If $(E, \mathcal{D})$ is an $A_\infty$-extension of $B$ by $A$, then $\mathcal{D} \circ \Xi$ defines an OCHA structure 
on the pair $(B,A)$. The $L_\infty$-structure on $B$ is the Lada-Markl symmetrization of the $A_\infty$-structure on $B$.
\end{thm*}
\noindent If the $A_\infty$-extension $(E, \mathcal{D})$ splits, then the OCHA structure induced by the above theorem 
reduces to an $A_\infty$-algebra over an $L_\infty$-algebra. Such structures were introduced by Kajiura and Stasheff
in \cite{KS06a} as the strong homotopy version of actions by derivations of Lie algebras on associative algebras.
Structures containing pairs $(L,A)$ where the Lie algebra $L$ acts by derivations on the associative algebra $A$ have 
appeared in different contexts in the literature, see \cite{Hue04} for an overview. 
If $A$ is commutative, $L$ is an $A$-module and $L$ acts by derivations
on $A$, then the pair $(L,A)$ is called a {\it Lie-Rinehart Algebra}, see \cite{Hue90}. 
In the case where $A$ is not necessarily commutative and $L$ does not need to be an $A$-module,
the pair is called a {\it Leibniz Pair} by Flato, Gerstenhaber and Voronov \cite{FGV95}. Thus, 
$A_\infty$-algebras over $L_\infty$-algebras are {\it Strong Homotopy Leibniz Pairs}.

The universal enveloping $A_\infty$-algebra of an OCHA is introduced in Section \ref{universal}. 
Its definition is completely analogous to the universal enveloping $A_\infty$-algebra of 
an $L_\infty$-algebra introduced by Lada and Markl in \cite{LM95}.
Given an OCHA 
$(\cls,\opn, \mathcal{D})$, 
its universal enveloping $A_\infty$-algebra is denoted $\mathcal{U}_\infty(\cls,\opn)$, while 
$\mathcal{U}_\infty(\cls)$ denotes Lada-Markl $A_\infty$-algebra. 
In section \ref{universal} we prove the following theorem:
\begin{thm*}
The universal enveloping $A_\infty$-algebra $\mathcal{U}_\infty(\cls,\opn)$ of an OCHA $(\cls, \opn, \mathcal{D})$ 
is an $A_\infty$-extension of $\mathcal{U}_\infty(\cls)$ by $\langle \opn \rangle$: 
\[ 0 \to \langle \opn \rangle \to \mathcal{U}_\infty(\cls,\opn) \to \mathcal{U}_\infty(\cls) \to 0 \]
where $\langle \opn \rangle$ denotes the $A_\infty$-ideal generated by $\opn$.
\end{thm*}
\noindent We close Section \ref{coshuff} by showing that the universal enveloping 
$A_\infty$-algebra of an OCHA satisfies a universal property naturally described in terms of $A_\infty$-extensions.

%\begin{thm*}
%Let $(\cls, \opn)$ be an OCHA. There is an $A_\infty$-extension $\mathcal{U}_\infty(\cls, \opn)$ of 
%$\mathcal{U}_\infty(\cls)$ by $\opn$  
% such that for any $A_\infty$-algebras $A$, $B$ and any $A_\infty$-extension $E$ of 
%$B$ by $A$ and any linear map $\cls \oplus \opn \stackrel{f_c \oplus f_o}{\longrightarrow} B \oplus A$ there exists 
%a unique morphism of $A_\infty$-extensions $\mathcal{U}_\infty(\cls, \opn) \to E$ such that the following diagram is commutative: 
%\begin{displaymath}
%\xymatrix{
%                                                         & \mathcal{U}_\infty(\cls, \opn) \ar[d] \\
% \cls \oplus \opn \ar[ru]^{\iota}\ar[r]^{\quad f_c \oplus f_o} &      E                                   }
%\end{displaymath}
%where $\iota : \cls \oplus \opn \to \mathcal{U}_\infty(\cls, \opn)$ is the inclusion.
%\end{thm*}

\subsection*{Notation and Conventions}
Let us fix a field $k$ of characteristic zero. By a vector space we will always mean a $\mathbb{Z}$-graded vector space.
Let $V$ be a vector space, we define a left action of the symmetric group $S_n$ on $V^{\otimes n}$ in the 
following way:  if $\tau \in S_2$ is a transposition, then the action is given by 
$x_1 \otimes x_2 \stackrel{\tau}{\mapsto} (-1)^{|x_1||x_2|} x_2 \otimes x_1$. Since any $\sigma \in S_n$ is a 
composition of transpositions, the sign of the action of $\sigma$ on $V^{\otimes n}$ is well defined: 
\begin{equation*}
 x_1 \otimes \cdots \otimes x_n \stackrel{\sigma}{\mapsto} (-1)^{\epsilon(\sigma)}
x_{\sigma(1)} \otimes \cdots \otimes x_{\sigma(n)}. 
\end{equation*}
We will refer to $(-1)^{\epsilon(\sigma)}$ as the Koszul sign of the permutation. Let us 
define $(-1)^{\chi(\sigma)} = (-1)^\sigma (-1)^{\epsilon(\sigma)}$, where $(-1)^\sigma$ is the 
sign of the permutation. 
Given two homogeneous maps $f,g : V \rar W$,
we will follow the Koszul sign convention for the tensor product:
$(f \otimes g)(v_1 \otimes v_2) = (-1)^{|g||v_1|} (f(v_1) \otimes g(v_2))$.

In general, a Strong Homotopy Algebra defined on some vector space $V$ will be denoted as 
a pair $(V, \mathcal{D})$. The symbol $\mathcal{D}$ stands for the SH-structure. In this work, $\mathcal{D}$ can be thought 
of as a sequence of multilinear operations as well as a coderivation differential on some coalgebra.  
We observe that any vector space endowed with some SH structure has in particular a differential operator
which is part of the SH structure. 

\section{Open Closed Homotopy Algebras}
Here we define OCHA using the same grading and signs conventions of \cite{KS06a} which is 
the more appropriate for the coalgebra description. For a equivalent description where the 
grading and signs have a geometrical meaning, see \cite{Hoe06a}. 
Let us begin by recalling the definition of $L_\infty$-algebras \cite{LS93}.
\begin{defi}[Strong Homotopy Lie algebra]
A strong homotopy Lie algebra (or $L_\infty$-algebra) is a 
$\mathbb{Z}$-graded vector space $V$ endowed 
with a collection of degree one graded symmetric $n$-ary brackets $l_n : V^{\otimes n} \rar V$, 
such that $l_1^{^2} = 0$ and for $n \geqs 2$:
\begin{equation}\label{homotopy_lie}
  \sum_{\sigma \in \Sigma_{k+l=n} } \!\!\!\! \epsilon(\sigma)\; 
 l_{l}(l_k(v_{\sigma(1)}, \dots,v_{\sigma(k)}),v_{\sigma(k+1)}, \dots,v_{\sigma(n)}) = 0 
\end{equation} 
where $\sigma$ runs over all $(k,l)$-unshuffles, i.e., permutations $\sigma \in S_n$ 
such that $\sigma(i) < \sigma(j)$ for $1 \leqs i < j \leqs k$ and for  
$k+1 \leqs i < j \leqs k+l$. 
\end{defi}  

\begin{defi}[Open-Closed Homotopy Algebra $-$ OCHA] \label{OCHA_def}
An OCHA consists of a pair of graded vector spaces $(\cls,\opn)$ endowed with
two sequences of degree one multilinear operations 
$\mathfrak{l} = \{ l_n : L^{\otimes n} \rar L \}_{n \geqs 1}$ and 
$\mathfrak{n} = \{ n_{p,q} : L^{\otimes p} \otimes A^{\otimes q} \rar A \}_{p+q \geqs 1}$, 
such that $(L,\mathfrak{l})$ is an $L_\infty$-algebra
and the two families satisfy the following compatibility condition: 
\begin{multline*}
 \mquad \sum_{\scriptscriptstyle{\sigma \in \Sigma_{p+r=n}, \;p \geqs 2}}
 \mquad (-1)^{\epsilon(\sigma)}
n_{1+r,m}(l_p(v_{\sigma(1)}, \dots , v_{\sigma(p)}),v_{\sigma(p+1)}, \dots ,
                                          v_{\sigma(n)},a_1, \dots, a_m)\; + \\
\hspace*{-10.5em}
 + \mquad \; \sum_{\begin{array}{c}
              \scriptscriptstyle{\sigma \in \Sigma_{p+r=n}} 
              \\[-1ex] 
              \scriptscriptstyle{i+j=m-s}
            \end{array}} 
  (-1)^{\mu_{p,i}(\sigma)} n_{p,i+1+j}(v_{\sigma(1)}, . . , v_{\sigma(p)},a_1, . . , a_i, \\[-5ex]
 n_{r,s}(v_{\sigma(p+1)}, . . , v_{\sigma(n)},a_{i+1}, . . , a_{i+s}), a_{i+s+1}, . . , a_m) = 0.
\end{multline*}
\begin{multline*}
\mquad\mbox{ where } \mu_{p,i}(\sigma) = s + i + si + ms +\epsilon(\sigma) + 
s(v_{\sigma(1)} + \cdots + v_{\sigma(p)} + a_1 + \cdots + a_i) + \\ + 
(a_1 + \cdots + a_i)(v_{\sigma(i+1)}) + \cdots + v_{\sigma(n)}). 
\end{multline*}
\end{defi}

\section{Coderivations}\label{coderivations}
Here we shall briefly recall the notion of coderivation and how the Gerstenhaber 
bracket can be seen intrinsically as the graded commutators of coderivations on the tensor coalgebra,
see \cite{Sta93} for details. 

Given a coalgebra $(C, \Delta, \varepsilon)$, a coderivation of $C$ is a linear 
map $f: C \rar C$ such that $(f \otimes 1 + 1 \otimes f)\Delta = \Delta f$ and 
$\varepsilon f = 0$. 
In the case of the tensor coalgebra $T^c V$ cogenerated by $V$, any linear 
map $f : V^{\otimes k} \rar V$ can be lifted to a coderivation 
$\hat{f} : T^c V \rar T^c V$ defined as $\hat{f}(v_1 \otimes \cdots \otimes v_n) = 0$
if $n < k$ and 
%\[ \hat{f}(v_1 \otimes \cdots \otimes v_n) = 
%\sum_{i=0}^{n-k} v_1 \otimes \cdots \otimes v_i \otimes 
%f(v_{i+1} \otimes \cdots \otimes v_{i+k}) \otimes v_{i+k+1} \otimes \cdots \otimes v_n \]  
\begin{equation}\label{lifting_tensor}
\hat{f}(v_1 \otimes \cdots \otimes v_n) = 
\sum_{i=0}^{n-k} (1^{\otimes i} \otimes f \otimes 1^{n-i-k})(v_1 \otimes \cdots \otimes v_n), 
\quad {\rm for } \; n \geqs k.
\end{equation} 

Consequently, any map $f : T^c V \rar V$ can be lifted to $\hat{f} : T^c V \rar T^c V$
by adding up the liftings of each component of the map $f$. The lifting as a coderivation 
defines an isomorphism between the vector spaces ${\rm Hom}(T^cV,V)$ and 
${\rm Coder}(T^cV)$. 

Restricting attention to graded coderivations, we have a space with the structure 
of a graded Lie algebra with bracket given by the graded commutator of compositions: 
$[\theta, \phi] = \theta \circ \phi - (-1)^{|\theta||\phi|} \phi \circ \theta$. 
Under the isomorphism provided by the lifting as a coderivation, the above bracket induces
a graded Lie algebra structure on $\bigoplus_{n \geqs 0}{\rm Hom}(V^{\otimes n},V)$ (the 
graded space corresponding to the space of graded coderivations). Stasheff has shown that 
the bracket induced through the above isomorphism is, up to sign, the Gerstenhaber bracket \cite{Sta93}. 
The sign can be adjusted using suspension and desuspension operators, see \cite{Sta93,DMZ07}.

\begin{obs}
It is also possible to define 
the lifting as a coderivation for the symmetric coalgebra $\Lambda^c V$ as well as for 
$\Lambda^c V \otimes T^c U$, see formulas on page \pageref{tensor}. However,
those formulas holds only in the particular case of the symmetric and tensor coalgebras.
For a more general study of coderivations, one can consider coderivations of a $\mathcal{P}$-coalgebra $X$, where $\mathcal{P}$ is an operad \cite{MSS02}. 
\end{obs}

\begin{obs}
  Since the spaces $\cls$ and $\opn$ are graded, there are two ways of grading the space 
  of coderivations. One way is through the \emph{external} degree, i.e., the grading 
  by the number of variables. The second is induced by the grading
  of $\cls \otimes \opn$ and is called \emph{internal} degree. In this paper we will 
  only use the internal degree. 
\end{obs}

\subsection{Coderivations on $\Lambda^c \mathcal{H}_c \otimes T^c \mathcal{H}_o$}\label{coder_ocha}
According to its operadic description \cite{KS06a, KS06b, Hoe06a}, each OCHA operation $n_{p,q}$ and $l_n$ correspond to 
a partially planar corolla: 
\begin{equation}\label{corollae}
 \raisebox{-3em}{\begin{picture}(0,0)%
\includegraphics{lk.pstex}%
\end{picture}%
\setlength{\unitlength}{3947sp}%
\begingroup\makeatletter\ifx\SetFigFont\undefined%
\gdef\SetFigFont#1#2#3#4#5{%
  \reset@font\fontsize{#1}{#2pt}%
  \fontfamily{#3}\fontseries{#4}\fontshape{#5}%
  \selectfont}%
\fi\endgroup%
\begin{picture}(1814,1072)(8938,-829)
\put(9382, 95){\makebox(0,0)[lb]{\smash{{\SetFigFont{8}{14.4}{\rmdefault}{\mddefault}{\updefault}$1$}}}}
\put(9790, 80){\makebox(0,0)[lb]{\smash{{\SetFigFont{8}{14.4}{\rmdefault}{\mddefault}{\updefault}$2$}}}}
\put(9982,-98){\makebox(0,0)[lb]{\smash{{\SetFigFont{8}{14.4}{\rmdefault}{\mddefault}{\updefault}$\dots$}}}}
\put(10055,111){\makebox(0,0)[lb]{\smash{{\SetFigFont{8}{14.4}{\rmdefault}{\mddefault}{\updefault}$\dots$}}}}
\put(10476, 81){\makebox(0,0)[lb]{\smash{{\SetFigFont{8}{14.4}{\rmdefault}{\mddefault}{\updefault}$n$}}}}
\put(9200,-444){\makebox(0,0)[lb]{\smash{{\SetFigFont{8}{14.4}{\rmdefault}{\mddefault}{\updefault}$l_n = $}}}}
\end{picture}%
 \quad \begin{picture}(0,0)%
\includegraphics{nkl.pstex}%
\end{picture}%
\setlength{\unitlength}{3947sp}%
\begingroup\makeatletter\ifx\SetFigFont\undefined%
\gdef\SetFigFont#1#2#3#4#5{%
  \reset@font\fontsize{#1}{#2pt}%
  \fontfamily{#3}\fontseries{#4}\fontshape{#5}%
  \selectfont}%
\fi\endgroup%
\begin{picture}(2035,1148)(6498,-3097)
\put(7269,-2109){\makebox(0,0)[lb]{\smash{{\SetFigFont{8}{14.4}{\rmdefault}{\mddefault}{\updefault}$\dots$}}}}
\put(8017,-2096){\makebox(0,0)[lb]{\smash{{\SetFigFont{8}{14.4}{\rmdefault}{\mddefault}{\updefault}$\dots$}}}}
\put(7113,-2141){\makebox(0,0)[lb]{\smash{{\SetFigFont{8}{14.4}{\rmdefault}{\mddefault}{\updefault}$1$}}}}
\put(7870,-2144){\makebox(0,0)[lb]{\smash{{\SetFigFont{8}{14.4}{\rmdefault}{\mddefault}{\updefault}$1$}}}}
\put(6800,-2626){\makebox(0,0)[lb]{\smash{{\SetFigFont{8}{14.4}{\rmdefault}{\mddefault}{\updefault}$ n_{p,q} = $}}}}
\put(7246,-2266){\makebox(0,0)[lb]{\smash{{\SetFigFont{8}{14.4}{\rmdefault}{\mddefault}{\updefault}$\dots$}}}}
\put(7885,-2266){\makebox(0,0)[lb]{\smash{{\SetFigFont{8}{14.4}{\rmdefault}{\mddefault}{\updefault}$\dots$}}}}
\put(8257,-2105){\makebox(0,0)[lb]{\smash{{\SetFigFont{8}{14.4}{\rmdefault}{\mddefault}{\updefault}$q$}}}}
\put(7523,-2116){\makebox(0,0)[lb]{\smash{{\SetFigFont{8}{14.4}{\rmdefault}{\mddefault}{\updefault}$p$}}}}
\end{picture}%
.} 
\end{equation}
In the language of trees, the OCHA constraint can be stated as follows: 
``The OCHA operad has no corolla with spatial root and planar leaves''.
The Axelrod-Singer compactification of the moduli space of points on the closed disc is a manifold with corners 
whose boundary is labelled by the partially planar trees obtained after grafting corollae of the above type. 
Other types of partially planar corollae will not appear in that boundary strata. That is the geometrical origin 
of the OCHA constraint.  

In this section we show that any coderivation on $\Lambda^c \mathcal{H}_c \otimes T^c\mathcal{H}_o$ 
satisfies the OCHA constraint, i.e., that the OCHA structure is {\it intrinsic} to the coalgebra 
$\Lambda^c \mathcal{H}_c \otimes T^c\mathcal{H}_o$. 
%n other words, we show that the constraint on the OCHA structure maps $\{ n_{p,q} \}$ and $\{ l_n \}$, 
%whose geometric meaning we have discussed in section, are in fact intrinsic to the above coalegebra.
%%We will use a natural injective coalgebra map 
%%$\Xi : \Lambda^c U \otimes T^c V \rar T^c(U \oplus V)$ and formulas 
%%(\ref{lift1}) and (\ref{lift2}) on the appendix. The definition 
%%of $\Xi$ involves the shuffle product.
We begin by recalling the definition of the shuffle product.
\begin{defi}
For any vector space $E,$ let 
$a_1 \otimes \cdots \otimes a_n \in E^{\otimes n}$ and $0 \leqs i \leqs n$. 
The \emph{shuffle product} of 
$a_1 \otimes \cdots \otimes a_i \in E^{\otimes i}$ and 
$a_{i+1} \otimes \cdots \otimes a_n \in E^{\otimes n-i}$
is:
\begin{equation}
 Sh_{i,n-i}(a_1 \otimes \cdots \otimes a_i|a_{i+1} \otimes \cdots \otimes a_n) := \mquad
\sum_{\sigma \in \mathfrak{U}_{(i,n-i)}} (-1)^{\epsilon(\sigma)} a_{\sigma(1)} \otimes \dots \otimes a_{\sigma(n)}
\end{equation} 
where $\mathfrak{U}_{(i,n-i)}$ is the set of all $(i,n-i)$-shuffles, i.e., 
$\sigma \in \Sigma_n$ with
\[ \sigma^{-1}(1) < \cdots < \sigma^{-1}(i), \ \sigma^{-1}(i+1) < \cdots < \sigma^{-1}(n). \]   
\end{defi}
The shuffle product $Sh : T^cE \otimes T^cE \to T^cE$ defines an associative product that is compatible 
with the deconcatenation coproduct $\Delta : T^cE \to T^cE \otimes T^cE$, so $Sh$ is a coalgebra map with 
respect to $\Delta$. 

Now consider $E = \mathcal{H}_c \oplus \mathcal{H}_o$ and define 
$\Xi: \Lambda^c \cls \otimes T^c \opn \rar T^c(\cls \oplus \opn):$ 
\[ \Xi((c_1 \wedge \cdots \wedge c_p) \otimes (o_1 \otimes \cdots \otimes o_q)) = 
Sh(\chi(c_1 \wedge \cdots \wedge c_p)|(o_1 \otimes \cdots \otimes o_q)) \]
where $\chi : \Lambda^c \cls \rar T^c \cls$ is the \emph{symmetrization coalgebra map}: 
\begin{equation*} \label{chi}
 \chi(c_1 \wedge \cdots \wedge c_n) = \sum_{\sigma \in S_n} (-1)^{\epsilon(\sigma)}
                                      c_{\sigma(1)} \otimes \cdots \otimes c_{\sigma(n)}.
\end{equation*} 
 
We will sometimes refer to $\Xi: \Lambda^c \cls \otimes T^c \opn \rar T^c(\cls \oplus \opn):$ 
as the {\it symmetrization and shuffling coalgebra map}.

Since we are working over a field of characteristic zero, the symmetrization and shuffling $\Xi$ is an 
injective coalgebra map. Now let 
$p_n : \Lambda^c \cls \otimes T^c \opn \rar \bigoplus_{p+q=n}\cls^{\wedge p} \otimes \opn^{\otimes q}$ be 
the canonical projection. We observe
that $p_n = \Xi^{-1} \pi_n \Xi$ where $\pi_n : T^c (\cls \oplus \opn) \rar (\cls \oplus \opn)^{\otimes n}$ is the canonical projection. In particular, $p_1 = \pi_1 \Xi : \Lambda^c \cls \otimes T^c \opn \rar \cls \oplus \opn$. 
 It is not diFficult 
to see that $\pi_n = (\pi_1 \otimes \cdots \otimes \pi_1) \Delta_{\otimes}^{(n-1)}$, 
Where $\Delta_{\otimes}$ denotes the deconcatenation coproduct on $T^c (\cls \oplus \opn)$,
we thus have:
\begin{equation} \label{f+g}
\begin{split}
 p_n & = \Xi^{-1} (\pi_1 \otimes \cdots \otimes \pi_1) \Delta_{\otimes}^{(n-1)} \Xi \\
          & = \Xi^{-1} (\pi_1 \otimes \cdots \otimes \pi_1) \Xi^{\otimes n} 
                                                            \Delta^{(n-1)} \\
          & = \Xi^{-1} (p_1 \otimes \cdots \otimes p_1)\Delta^{(n-1)} 
\end{split}
\end{equation} 

\noindent where $\Delta$ is the 
coproduct on $\Lambda^c \cls \otimes T^c \opn$ defined by (\ref{lambdat}) 
in the appendix. 

Given a coderivation $\phi \in {\rm Coder}(\Lambda^c \cls \otimes T^c \opn)$, applying it 
to (\ref{f+g}):
\begin{equation*} 
\begin{split}
 p_n \phi 
   & = \Xi^{-1} (p_1 \otimes \cdots \otimes p_1)\Delta^{(n-1)} \phi \\
   & = \Xi^{-1} \sum (p_1 \otimes \cdots \otimes p_1 \phi
                     \otimes \cdots \otimes p_1)\Delta^{(n-1)},
\end{split}
\end{equation*}  
we conclude that $\phi$ is determined by its projection 
$p_1 \phi : \Lambda^c \cls \otimes T^c \opn \rar \cls \oplus \opn$, $p_1 \phi = g \oplus f$. 
We can write $g$ and $f$ as: $g = \sum g_{p,q}$ and $f = \sum f_{p,q}$, where 
$g_{p,q} : \cls^{\wedge p} \otimes \opn^{\otimes q} \rar \cls$ and 
$f_{p,q} : \cls^{\wedge p} \otimes \opn^{\otimes q} \rar \opn$. 

Let us show that $\phi$ 
is given by the lifting of $g_{p,q}$ and $f_{p,q}$ as coderivations:
\begin{equation*} 
\begin{split}
 p_n \phi 
   & = \Xi^{-1} \sum (p_1 \otimes \cdots \otimes (f + g)
                     \otimes \cdots \otimes p_1)\Delta^{(n-1)} \\
   & = \Xi^{-1} \sum (p_1 \otimes \cdots \otimes 
(\sum g_{p,q} + \sum f_{p,q}) \otimes \cdots \otimes p_1)\Delta^{(n-1)}, 
\end{split}
\end{equation*}  
consequently, $p_n \phi$ is the sum of the two expressions:
\begin{gather}
 \sum_{p,q}\Xi^{-1} \sum (p_1 \otimes \cdots \otimes \label{g_pq}
 g_{p,q} \otimes \cdots \otimes p_1)\Delta^{(n-1)} \\
 \sum_{p,q}\Xi^{-1} \sum (p_1 \otimes \cdots \otimes \label{f_pq}
 f_{p,q} \otimes \cdots \otimes p_1)\Delta^{(n-1)}. 
\end{gather}
Extending $f_{p,q} \mapsto \hat{f}_{p,q}$ by  formula (\ref{lift1}),
we see that (\ref{f_pq}) is equal to $p_n \sum_{p,q} \hat{f}_{p,q}$.  
On the other hand, we lift 
$g_{p,q} : \cls^{\wedge p} \otimes \opn^{\wedge q} \to \cls$ to
$\hat{g}_{p,q} : \Lambda^c \cls \otimes T^c \opn \rar \Lambda^c \cls \otimes T^c \opn$, 
with $\hat{g}_{p,q}$ given by:
\[ \hat{g}_{p,q}((u_1 \wedge \cdots \wedge u_n) \otimes 
 (v_1 \otimes \cdots \otimes v_m)) = 0, \qquad \mbox{if $n < p$ or $m < q$}, \] 
and
\begin{multline} \label{lift3}
 \hat{g}_{p,q}((u_1 \wedge \cdots \wedge u_n) \otimes (v_1 \otimes \cdots \otimes v_m)) = \\
 \sum_{\sigma \in \mathfrak{S}_{p,n-p}} \pm 
 (g_{p,q}(u_{\sigma(1)}, \dots, u_{\sigma(p)},v_{i+1}, \dots, v_{i+q}) 
 \wedge u_{\sigma(p+1)} \wedge \cdots \wedge u_{\sigma(n)}) 
 \otimes \hspace{\stretch{1}} \\[-1.3em]
 \otimes(v_1 \otimes \cdots \otimes \widehat{v_{i+1}} \otimes \cdots \otimes 
                                    \widehat{v_{i+q}} \otimes \cdots \otimes v_m)
\end{multline}
where 
$\pm = (-1)^{\displaystyle{ (\epsilon(\sigma) + (u_{\sigma(p+1)} + \cdots + u_{\sigma(n)} + v_{1} 
+ \dots + v_{i})(v_{i+1} + \dots + v_{i+q}))}}$ and $\hat{v}$ means that $v$ is omitted in 
the expression. 

 It follows that $p_n \sum_{p,q} \hat{g}_{p,q}$ is equal to expression
(\ref{g_pq}),  and the coderivation can thus be written as 
%\begin{equation}
$\phi = \sum \hat{g}_{p,q} + \sum \hat{f}_{p,q} \ .$ 
%\end{equation}
Observe that (\ref{lift3}) reduces to (\ref{lift2}) in the appendix when $q=0$. 

%In order to complete our characterization of all coderivations on    
%$\Lambda^c\cls \otimes T^c\opn$, we need to prove that, since $\phi$ is 
%a coderivation, the maps $g_{p,q}$ with $q \geqs 1$ have to be zero, i.e., 
%the coderivations are in OCHA form. 

%We have shown that any coderivation in Coder$(\Lambda^c\cls \otimes T^c\opn)$ has the form
%$\phi = \sum \hat{g}_{p,q} + \sum \hat{f}_{p,q}$. 

We say that a coderivation $\phi$ satisfies the {\it OCHA constraint} if
\begin{equation}
 \phi = \sum \hat{g}_p + \sum \hat{f}_{p,q} 
\end{equation}
for $g_p : \cls^{\wedge p} \to \cls$ and $f_{p,q} : \cls^{\wedge p} \otimes \opn^{\otimes q} \to \opn$.

\begin{prop} \label{coder}
Every coderivation $\phi \in {\rm Coder}(\Lambda^c \cls \otimes T^c \opn)$ satisfies the OCHA constraint.
\end{prop}

%\begin{prop} \label{coder}
% Every coderivation $\phi \in {\rm Coder}(\Lambda^c \cls \otimes T^c \opn)$ can be written as 	
% \[ \phi = \sum_{n \geqs 0} \hat{g}_n + \sum_{p,q \geqs 0} \hat{f}_{p,q} \] 
% where $\hat{g}_n$ and  $\hat{f}_{p,q}$
% are obtained by lifting, as defined by expressions (\ref{lift1}) and (\ref{lift2}), 
% multilinear maps $g_n : \cls^{\wedge n} \rar \cls$ and $f_{p,q} : \cls^{\wedge p} \otimes \opn^{\otimes q} \rar \opn$. 
%\end{prop}

\begin{proof}
We have already seen that formula (\ref{f+g}) implies that any coderivation 
$\phi \in {\rm Coder}(\Lambda^c \cls \otimes T^c \opn)$ is uniquely written as 
$\phi = \sum \hat{g}_{p,q} + \sum \hat{f}_{p,q}$ where 
$\hat{f}_{p,q}$ was defined in formula (\ref{lift1}) and $\hat{g}_{p,q}$ in formula (\ref{lift3}). 
We know that $\hat{f}_{p,q}$ is a coderivation for any $p,q$ and that $\hat{g}_{p,q}$ is a coderivation
if $q=0$. Thus $\phi - (\sum \hat{f}_{p,q} + \sum \hat{g}_{p,0})$ = \mquad
\raisebox{-1ex}{$\begin{array}{c} \sum  \\[-0.5ex] \scriptscriptstyle{(p,q) 
\geqs (0,1)} \end{array}$} \mquad
$\hat{g}_{p,q}$
is a coderivation:
\begin{equation} \label{proofcoder1}
 \sum_{\scriptscriptstyle (p,q) \geqs (0,1)} \!\!\!
 (\hat{g}_{p,q} \otimes 1 + 1 \otimes \hat{g}_{p,q})\bar{\Delta} = \mquad
 \sum_{\scriptscriptstyle (p,q) \geqs (0,1)} \!\!\!
 \bar{\Delta} \hat{g}_{p,q}
\end{equation}
where $\bar{\Delta}$ denotes the reduced comultiplication: $\Delta x = x \otimes 1 + \bar{\Delta} x + 1 \otimes x$. 

We will prove that $g_{p,q} \equiv 0$ for any $p \geqs 0, q \geqs 1$. 
To simplify the exposition we will use the notation: 
$u_{[p]} = u_1 \land \cdots \land u_p$. 
From (\ref{lift3}), we have
\begin{equation*} \mquad
 \hat{g}_{p,q}(u_{[p]} \otimes (v_1 \otimes \cdots \otimes v_{q+1})) = 
 g_{p,q}(u_{[p]}; v_1 \otimes \cdots \otimes v_{q}) \otimes v_{q + 1} \,\pm\,   
 g_{p,q}(u_{[p]}; v_2 \otimes \cdots \otimes v_{q+1}) \otimes v_{1}
\end{equation*}
where $\pm = (-1)^{|v_1|(|v_2| + \cdots + |v_{q+1}|)}$.

Applying $u_{[p]} \otimes (v_1 \otimes \cdots \otimes v_q \otimes v_{q+1})$ to both sides 
of equation (\ref{proofcoder1}) and projecting the result onto 
$(\cls \otimes \opn) \oplus (\opn \otimes \cls)$ (viewed as a subspace of 
$(\Lambda^c \cls \otimes T^c \opn) \otimes (\Lambda^c \cls \otimes T^c \opn)$) we have:
\begin{multline*}
g_{p,q}(u_{[p]}; v_1, \dots, v_q) \otimes v_{q+1} \pm
   v_1 \otimes g_{p,q}(u_{[p]}; v_2, \dots, v_{q+1}) = \\
\mquad = g_{p,q}(u_{[p]}; v_1, \dots, v_q) \otimes v_{q+1} \pm
   g_{p,q}(u_{[p]}; v_2, \dots, v_{q+1}) \otimes v_1 \,\pm \\
   \pm v_{q+1} \otimes g_{p,q}(u_{[p]}; v_1, \dots, v_q) \pm
   v_1 \otimes g_{p,q}(u_{[p]}; v_2, \dots, v_{q+1}).
\end{multline*}
where $\pm$ is, in order of appearance: $(-1)^{|v_1|(|u_{[p]}| + 1)}$; 
$(-1)^{|v_{q+1}|(|u_{[p]}| + |v_1| + \cdots + |v_q| + 1)}$;
$(-1)^{|v_1|(|v_2| + \cdots + |v_{q+1}|)}$; $(-1)^{|v_1|(|u_{[p]}| + 1|)}.$ 

\noindent
It follows that:
\[
g_{p,q}(u_{[p]}; v_2, \dots, v_{q+1}) \otimes v_1 
   \pm v_{q+1} \otimes g_{p,q}(u_{[p]}; v_1, \dots, v_q) = 0,
\]
one summand belongs to $(\cls \otimes \opn)$ while the other one belongs to $(\opn \otimes \cls)$, 
hence both are zero. Assuming $v_{q+1} \neq 0$, we have $g_{p,q}(u_{[p]}; v_1, \dots, v_q) = 0$.
%
%
%%\[ g_{p,q}(u_{[p]}; v_2, \dots, v_{q+1}) \otimes v_1 \; 
%%= \; v_{q+1} \otimes g_{p,q}(u_{[p]}; v_1, \dots, v_q) =  0 \]
%%for any $u_1 \wedge \cdots \wedge u_p \in \cls^{\wedge p}$ and 
%%$v_1 \otimes \cdots \otimes v_{q+1} \in \opn^{\otimes q+1}$. 
%
%We thus conclude that $g_{p,q}$ is identically zero for $p \geqs 0, q \geqs 1$. 
\end{proof}
\begin{thm}
  An OCHA structure on the pair $(\mathcal{H}_c, \mathcal{H}_o)$ is equivalent to a degree one 
  coderivation $D \in {\rm Coder}^1(\Lambda^c \mathcal{H}_c \otimes T^c \mathcal{H}_o)$ such that $D^2=0$.
\end{thm}

\subsection{OCHA-morphisms}

%We will now recall the notion of morphism between OCHAs.
Consider two families of degree zero linear maps $f_k : \Lambda^k \cls \to \cls'$, for $k \geqs 1$
and $f_{p,q} : \Lambda^p \cls \otimes \opn^{\otimes q} \to \opn'$ for $p+q \geqs 1$.
Given two OCHAs $(\cls \oplus \opn, D)$ and $(\cls' \oplus \opn', D')$, according Kajiura and Stasheff
we say that the maps $\{ f_k \}_{k \geqs 1}$ and $\{ f_{p,q} \}_{p+q \geqs 1}$ define 
an OCHA-morphism when they commute with the OCHA structures after lifted as a coalgebra morphism.
More precisely: 
\[ \frak{f} \circ D = D' \circ \frak{f} \]
where $\frak{f} : \Lambda^c \cls \otimes T^c\opn \to \Lambda^c \cls' \otimes T^c\opn'$ is the coalgebra 
morphism obtained by lifting the degree zero linear maps $\{ f_k \}_{k \geqs 1}$
and $\{ f_{p,q} \}_{p+q \geqs 1}$. Explicit formulas for OCHA-morphisms are 
available in \cite{KS06a}. In the particular case 
of linear OCHA-morphisms, explicit formulas  are provided below.

We say that an OCHA-morphism $\frak{f}$ is {\it linear} when it is obtained by lifting maps $g : \cls \to \cls'$, 
$f_{0,1} : \opn \to \opn'$ and $f_{1,0} : \cls \to \opn'$. Denoting the OCHA structures
by $D = \frak{l} + \frak{n}$ and $D' = \frak{l}' + \frak{n}'$, 
%and the components 
%of the linear OCHA morphism by $g$, $f_{0,1}$, $f_{1,0}$, 
equation $\frak{f} \circ D = D' \circ \frak{f}$ can be 
written as $g \circ l_n = l'_n \circ g^{\otimes n}$ for the $L_\infty$-structure maps and, for the 
remaining OCHA-structure maps as:
\begin{multline}
 f_{0,1}(\eta_{n,0}(c_1, \dots, c_n)) + f_{1,0}(l_n(c_1, \dots, c_n)) = \\ =
 \sum_{p+q=n} \frac{1}{p!} \  \eta'_{p,q}(g^{\otimes p} \otimes f_{1,0}^{\otimes q})\chi(c_1, \dots, c_n)  
\end{multline}
\begin{multline}
 f_{0,1}(\eta_{p,q}(c_1, \dots, c_p, o_1, \dots, o_q)) = \\ =\!\! \sum_{0 \leqs m \leqs p}\!\!  \frac{1}{m!} \ \eta'_{m,n}  
 (g^{\otimes m} \otimes Sh_{p-m,q}(f_{1,0}^{\otimes (p-m)} \otimes f_{0,1}^{\otimes q}) ) (\chi(c_1,.., c_p), o_1,.., o_q).
\end{multline}
The commutator $\chi$ and the shuffle product $Sh$ were defined in section \ref{coder_ocha}.
%where $i_k, j_k = 0,1$ are such that $m + i_1 + \cdots + i_n = p \mbox{ and } j_1 + \cdots + j_n = q$. 

\section{Commutators and shuffles of $A_\infty$-extensions} \label{coshuff}

Given an associative algebra with product $a \cdot b$, one can obtain a Lie 
algebra through the commutator:
$[a,b] = a \cdot b - b \cdot a$. This fundamental fact is also true in the context
of strongly homotopy Lie algebras. In fact, Lada and Markl \cite{LM95}
have used the symmetrization coalgebra map to relate 
$A_\infty$ and $L_\infty$ algebras and define the universal enveloping $A_\infty$-algebra of 
an $L_\infty$-algebra.

In this section we will show that analogous relations hold in the context of OCHA. 
In other words, we show that
an OCHA can be obtained from commutators and shuffles
of OCHA constrained 
$A_\infty$ structures on $\cls \oplus \opn$ or, equivalently, $A_\infty$-extensions
of $\cls$ by $\opn$. The universal enveloping $A_\infty$-algebra of an OCHA is then defined as an 
$A_\infty$-extension satisfying a natural universal property.

\subsection{Constraints on $A_\infty$-algebras as $A_\infty$-extensions}
In this section we show that the OCHA constraint on an $A_\infty$-algebra can be understood as 
an $A_\infty$-extension.

\begin{defi}
Let $(E,\mathcal{D})$ be an $A_\infty$-algebra which split in the category of vector spaces as: $E = A \oplus B$. 
We say that the $A_\infty$-algebra $E$ 
satisfy the OCHA constraint with respect to the splitting $E = A \oplus B$ if $\mathcal{D}$ has only components 
of the form: $\mathcal{D}^B_q : B^{\otimes q} \rightarrow B$ and  
$\mathcal{D}^A_{p,q} : A^{\otimes p} \otimes B^{\otimes q} \rightarrow A$.
\end{defi}

In this paper, $A_\infty$-extensions are defined using {\it linear} $A_\infty$-morphisms. 
A linear $A_\infty$-morphism between two $A_\infty$-algebras $(A, M=\{m_k\})$ and $(B,M'=\{m'_k\})$
is a degree zero linear map $f:A \to B$ such that:
\[ f \circ m_k = m'_k \circ (f \otimes \cdots \otimes f) \qquad \forall k \geqs 1. \]
In general, an $A_\infty$-morphism is given by a degree zero coalgebra morphism $\varphi : T^c(A) \to T^c(B)$ 
such that $M' \circ \varphi = \varphi \circ M$, viewing $M$ and $M'$ as coderivation differentials defining the 
$A_\infty$-structures.

\begin{defi}[$A_\infty$-ideal]
  Let $A$ be an $A_\infty$-algebra. An $A_\infty$-ideal of $A$ is a subespace $I \subseteq A$ such that, for any 
  $k \geqs 1$, $m_k(x_1, \dots, x_k) \in A$ whenever $x_i \in I$ for some $1 \leqs i \leqs k$.
\end{defi}

Notice that an $A_\infty$-ideal is in particular a subcomplex of $A$ and, 
if $f:A \to B$ is a linear $A_\infty$-morphism, then ${\rm Ker}(f) \subset A$ is 
an $A_\infty$-ideal of $A$. 

\begin{defi}
Let $A$ and $B$ be $A_\infty$-algebras. We say that an $A_\infty$-algebra $E$ is an $A_\infty$-extension of $B$ 
by $A$ if there exists an exact sequence 
\[ 0 \rightarrow A \rightarrow E \rightarrow B \rightarrow 0 \]
where each map is a linear $A_\infty$-morphism. 
\end{defi}

If $E$ is an extension of $B$ by $A$, then $E = A \oplus B$ as vector spaces and, since $A$ is an $A_\infty$-ideal in $E$,
we see that the $A_\infty$-algebra $E$ satisfies the OCHA constraint with respect to the splitting $E = A \oplus B$.  
Now that the appropriate concepts are given, we can just apply a well known argument \cite{Lada99,LM95} to prove the following theorem.
 
\begin{thm}\label{main}
  Let $A$ and $B$ be two $A_\infty$-algebras.
  If $(E, \mathcal{D})$ is an $A_\infty$-extension of $B$ by $A$, then $\mathcal{D} \circ \Xi$ defines an OCHA structure 
  on the pair $(B,A)$. The $L_\infty$ structure on $B$ is the Lada-Markl symmetrization of the $A_\infty$ structure on $B$.
\end{thm}
\begin{proof}
Since $\mathcal{D} = \{ \mathcal{D}_k \}$ is OCHA constrained, the composition 
$\mathcal{D} \circ \Xi = \{ \mathcal{D}_k \circ \Xi \}$ gives two sequence of maps: $l_k = \ell_k \circ \Xi : \Lambda^k \cls \rar \cls$ 
and $n_{p,q} = n_k \circ \Xi : \cls^{\wedge p} \otimes \opn^{\otimes q} \rar \opn$ with $p+q = k$.
Equations (\ref{lift1}) and (\ref{lift2}) tells us how to lift these maps as coderivations 
$\mathfrak{l}_k \in {\rm Coder}(\ocha)$ and $\mathfrak{n}_{p,q} \in {\rm Coder}(\ocha)$. 
We denote by $\widehat{\mathcal{D} \circ \Xi}$
the coderivation 
$\mathfrak{l} + \mathfrak{n} = \sum \mathfrak{l}_k + \sum \mathfrak{n}_{p,q} \in {\rm Coder}(\ocha)$.

Let us now consider the following diagram:
 \begin{gather}\label{Xidiagram}
  \xymatrix{ 
         \ocha \ar[r]^{\Xi} & T^c(\cls \oplus \opn)        &   \\
         \ocha \ar[u]^{\widehat{\mathcal{D} \circ \Xi}} \ar[r]^{\Xi} & 
                      T^c(\cls \oplus \opn) \ar[u]^{\hat{\mathcal{D}}} \ar[rd]^{\pi_c \oplus \pi_o} & \\
         \ocha \ar[u]^{\widehat{\mathcal{D} \circ \Xi}} \ar[r]^{\Xi} & 
                      T^c(\cls \oplus \opn) \ar[u]^{\hat{\mathcal{D}}} \ar[r]^{\mathcal{D}} & \cls \oplus \opn, }
         \\[-1.6em]\nonumber 
 \end{gather}
Since $\mathcal{D}$ is OCHA constrained, it is not dificult to check that 
$\Xi \circ \widehat{\mathcal{D} \circ \Xi} = \widehat{\mathcal{D}} \circ \Xi$, 
so the diagram is commuatative. Since $\mathcal{D} = \{ \mathcal{D}_k \}$ defines an $A_\infty$-algebra structure, 
we have ${\widehat{\mathcal{D}}}{}^{{}^{\scriptstyle 2}} = 0$. Since $\Xi$ is injective, from the above diagram 
we have $(\widehat{\mathcal{D} \circ \Xi})^2 = (\mathfrak{l} + \mathfrak{n})^2 = 0$.
\end{proof}

\begin{obs}\label{eoc}
If $E$ is an $A_\infty$-extension of $B$ by $A$, we will denote the OCHA obtained through commutators and shuffles
by $(E)_{OC}$.
\end{obs}

\subsection{Universal Enveloping $A_\infty$-algebra of an OCHA}\label{universal} 

Let us now construct the universal enveloping $A_\infty$ algebra of an OCHA.
Given an OCHA structure $\frak{l} + \frak{n} = \sum_{n \geqs 1} l_n + \sum_{p+q \geqs 1} n_{p,q} $
on a pair $(\cls, \opn)$, let $\mathcal{F}_\infty(\cls \oplus \opn)$ be the free $A_\infty$ algebra
generated by $\cls \oplus \opn$ with $A_\infty$ structure maps denoted by 
$\mu_n : (\cls \oplus \opn)^{\otimes n} \to \cls \oplus \opn$. %$\{ \mu_n \}_{n \geqs 1}$. 

We define the universal enveloping $A_\infty$-algebra $\mathcal{U}_\infty(\cls,\opn)$ as the quotient 
of $\mathcal{F}_\infty(\cls \oplus \opn)$ by the $A_\infty$-ideal $I$ generated by the relations:
\begin{gather}
 \!\!\!\! \sum_{\sigma \in S_p} (-1)^{\epsilon(\sigma)} \mu_{p+q}(Sh(c_{\sigma(1)}, .., c_{\sigma(p)}|o_1, .., o_q)) 
   = n_{p,q}(c_1, ..,c_p,o_1, .., o_q) \label{rel1} \\
 \!\!\!\! \sum_{\sigma \in S_p} (-1)^{\epsilon(\sigma)} \mu_{p}(c_{\sigma(1)}, \dots, c_{\sigma(p)}) 
   = l_p(c_1, \dots,c_p) + n_{p,0}(c_1, \dots,c_p), \label{rel2}
\end{gather}
for $p \geqs 0$, $q \geqs 1$ in (\ref{rel1}) and $p \geqs 1$ in (\ref{rel2}), where $c_i \in \cls$ and $o_j \in \opn$. 

In case $\opn = 0$, the OCHA structure reduces to an $L_\infty$-algebra structure on $\cls$ and the above 
construction reduces to the universal enveloping $A_\infty$-algebra $\mathcal{U}_\infty(\cls)$ of an 
$L_\infty$-algebra introduced by Lada and Markl in \cite{LM95}. In general, we have the following result 
relating the two constructions.
%Since
%$\frak{l} + \frak{n} = \sum_{n \geqs 1} l_n + \sum_{p+q \geqs 1} n_{p,q}$, we have: 
%\begin{align*}
%& (\frak{l} + \frak{n})(c_1, \dots,c_p,o_1, \dots, o_q) = n_{p,q}(c_1, \dots,c_p,o_1, \dots, o_q), \mbox{ if $q > 0$ } & \mbox{ and } & \\
%& (\frak{l} + \frak{n})(c_1, \dots,c_p,o_1, \dots, o_q) = l_p(c_1, \dots,c_p) + n_{p,0}(c_1, \dots,c_p), \mbox{ if $q=0$. } & \ &
%\end{align*}
\begin{thm}
The universal enveloping $A_\infty$-algebra $\mathcal{U}_\infty(\cls,\opn)$ 
of an OCHA $(\cls, \opn, \mathcal{D})$ 
is an $A_\infty$-extension of $\mathcal{U}_\infty(\cls)$ by $\langle \opn \rangle$: 
\[ 0 \to \langle \opn \rangle \to \mathcal{U}_\infty(\cls,\opn) \to \mathcal{U}_\infty(\cls) \to 0 \]
where $\langle \opn \rangle$ is the $A_\infty$-ideal generated by $\opn$.
\end{thm}
\begin{proof}
We just need to show that $\mathcal{U}_\infty(\cls,\opn)/\langle \opn \rangle$ is isomorphic to $\mathcal{U}_\infty(\cls)$, i.e., 
that $\mathcal{U}_\infty(\cls,\opn)/\langle \opn \rangle$ satisfies the universal property defining 
$\mathcal{U}_\infty(\cls)$.
Let $A$ be an $A_\infty$-algebra %with structure maps $\{ \hat\mu_n \}_{n \geqs 1}$ 
and let $f: \cls \to A$ be any linear map 
inducing a $L_\infty$-morphism from $\cls$ to $(A)_L$ (the $L_\infty$-algebra defined by commutators of $A$). 
Since $\mathcal{F}_\infty(\cls \oplus \opn)$ is free, there is a unique
$A_\infty$-morphism from $\mathcal{F}_\infty(\cls \oplus \opn)$ to $A$ extending $f$ and vanishing on $\langle \opn \rangle$. 
To see that it also vanishes on the ideal $I$ of relations defining 
$\mathcal{U}_\infty(\cls,\opn)$, we just need to note that it satisfies relations (\ref{rel1}) because it vanishes on 
$\langle \opn \rangle$ and relations (\ref{rel2}) since $f: \cls \to (A)_L$ is an $L_\infty$-morphism.
\end{proof}

The universal property characterizing $\mathcal{U}_\infty(\cls, \opn)$ is described as follows.
Let $(\cls, \opn, \mathcal{D})$ be an OCHA and 
%There is an $A_\infty$-extension $\mathcal{U}_\infty(\cls, \opn)$ of 
%$\mathcal{U}_\infty(\cls)$ by $\opn$  
let $A$ and $B$ be $A_\infty$-algebras. For any $A_\infty$-extension $E$ of 
$B$ by $A$ and any linear map $\cls \oplus \opn \stackrel{f}{\longrightarrow} E$ such 
that $\cls \oplus \opn \stackrel{f}{\longrightarrow} (E)_{OC}$ is a linear OCHA-morphism, there exists 
a unique morphism of $A_\infty$-extensions $\varphi : \mathcal{U}_\infty(\cls, \opn) \to E$ such that the following diagram is commutative:
\begin{displaymath}
\xymatrix{
                                            & \mathcal{U}_\infty(\cls, \opn) \ar[d]^{\varphi} \\
 \cls \oplus \opn \ar[ru]^{\iota}\ar[r]^{\quad f} &      E                                   }
\end{displaymath}
where $\iota : \cls \oplus \opn \to \mathcal{U}_\infty(\cls, \opn)$ is the inclusion.
Since $\mathcal{F}_\infty(\cls \oplus \opn)$ is free, there is a unique $A_\infty$-morphism 
$\varphi: \mathcal{F}_\infty(\cls \oplus \opn) \to E$ extending $f$. It vanishes on the ideal of OCHA relations 
(\ref{rel1}) and (\ref{rel2}) because $f: \cls \oplus \opn \to E_{OC}$ is a linear OCHA morphism. Hence $\varphi$ is 
well defined on $\mathcal{U}_\infty(\cls, \opn)$. This proves that there is a unique 
$A_\infty$-morphism $\varphi : \mathcal{U}_\infty(\cls, \opn) \to E$ extending $f$. 

To see that $\varphi$ is a morphism of $A_\infty$-extensions, we just need to observe that the following 
diagram is commutative:
\begin{displaymath}
\xymatrix{
  0 \ar[r] & \langle \opn \rangle \ar[d]%_{f_{0,1}}
  \ar[r]   & \mathcal{U}_\infty(\cls, \opn) \ar[r]\ar[d]_{\varphi} & \mathcal{U}_\infty(\cls)\ar[d]\ar[r] & 0 \\
  0 \ar[r] & A \ar[r] &  E \ar[r]                                             &   B        \ar[r]                    & 0.   }
\end{displaymath}
In fact, $\varphi : \mathcal{U}_\infty(\cls, \opn) \to E$ is an extension of an OCHA morphism $f : \cls \oplus \opn \to (E)_{OC}$. 
So, it respects the OCHA constraint taking the ideal $\langle \opn \rangle$ into the ideal $A$ and is thus well defined on the 
quotient $\mathcal{U}_\infty(\cls, \opn)/\langle \opn \rangle \to E/A$.

\subsection*{Acknowledgments}
The author wishes to thank Jim Stasheff and Murray Gerstenhaber for the kind hospitality
during his stay as a visiting graduate student 
at the University of Pennsylvania (CNPq-Brasil grant SWE-201064/04). 
We are also grateful to J. Stasheff, M. Gerstenhaber, H. Kajiura, T. Lada and M. Markl
for valuable discussions. 

\section*{Appendix: Lifting as a coderivation} \label{tensor}
Here we provide some formulas for the lifting as a coderivation in the case of the 
coalgebras $\Lambda^c U$ and $\Lambda^c U \otimes T^c \opn$.
$U$ and $\opn$ are $\mathbb{Z}$-graded vector spaces. The coproduct $\Delta$ on 
$\Lambda^c U \otimes T^c \opn$ is given explicitly by:
\begin{multline} \label{lambdat}
 \Delta ((u_1 \wedge \cdots \wedge u_m ) \otimes (v_1 \otimes \cdots \otimes v_n)) = \\
 = \sum_{\begin{array}{c}
        \scriptstyle 0 \leqs p \leqs m \\[-1ex]
        \scriptstyle 0 \leqs q \leqs n
       \end{array}}
 \sum_{\sigma \in \mathfrak{S}_{p,m-p}} (-1)^{\epsilon(\sigma)} (-1)^{\eta(p,q)}    
 ((u_{\sigma(1)} \wedge \cdots \wedge u_{\sigma(p)})  \otimes (v_1 \otimes \cdots \otimes v_q)) \otimes 
 \hspace{\stretch{1}} \\[-2.5em]
 \hspace{\stretch{1}}
 \otimes((u_{\sigma(p+1)} \wedge \cdots \wedge u_{\sigma(m)} ) \otimes (v_{q+1} \otimes \cdots \otimes v_n)), \\ 
\end{multline}
where :
$ \eta(p,q) = (u_{\sigma(p+1)} + \cdots + u_{\sigma(m)})(v_1 + \cdots + v_q).$

Given a map $f : {U}^{\wedge p} \otimes {\opn}^{\otimes q} \rar \opn$, we may lift it as a coderivation in 
the following way: for $r \geqs p,  s \geqs q$ we define: 
\begin{multline}\label{lift1}
 \hat{f}((u_1 \wedge \cdots \wedge u_r ) \otimes (v_1 \otimes \cdots \otimes v_s))  = \\
 = \sum_{\begin{array}{c}
        \scriptstyle \sigma \in \mathfrak{S}_{r-p,p} \\[-1ex]
        \scriptstyle 0 \leqs j \leqs s - q
       \end{array}} (-1)^{\mu_{r-p,j}(\sigma)}
 (u_{\sigma(1)} \wedge \cdots \wedge u_{\sigma(r-p)})  
 \otimes (v_1 \otimes \cdots \otimes v_j \otimes \hspace{\stretch{1}} \\[-2.5em]
 \hspace{\stretch{1}} \otimes f(u_{\sigma(r-p+1)}, \dots , u_{\sigma(r)},v_{j+1}, \dots, v_{j+q})\otimes \cdots \otimes v_s) 
\end{multline} 
\begin{multline*}
\hspace*{-.9em} \mbox{where: } 
\mu_{p,q}(\sigma) = \epsilon(\sigma) + (u_{\sigma(1)} + \cdots 
        + u_{\sigma(p)}) + (v_1 + \cdots + v_q) +  \\ 
        + (v_1 + \cdots + v_q)(u_{\sigma(q+1)} + \cdots + u_{\sigma(n)}).
\end{multline*}
It is not difficult to check that $\hat{f}$ is a coderivation.

Recall that a map $g: U^{\wedge p}  \rar U$  
may be lifted as a coderivation $\hat{g} : \Lambda^c U \rar \Lambda^c U$
so that: $\hat{g}(u_1 \wedge \cdots \wedge u_n) = 0$
for $n < p$ and for $n \geqs p$ is defined by:
\[ \hat{g}(u_1 \wedge \cdots \wedge u_n) = \sum_{\sigma \in \mathfrak{S}_{p,n-p}} (-1)^{\epsilon(\sigma)}
g(u_{\sigma(1)} \wedge \cdots \wedge u_{\sigma(p)}) \wedge u_{\sigma(p+1)} \wedge \cdots \wedge u_{\sigma(n)}. \] 
$g$ can be lifted as a coderivation of $\Lambda^c U \otimes T^c \opn$ by tensoring 
the above map with the identity of $T^c\opn$. We thus have       
$\hat{g} : \Lambda^c U \otimes T^c \opn \rar \Lambda^c U \otimes T^c \opn$
\begin{multline}\label{lift2}
 \hat{g}((u_1 \wedge \cdots \wedge u_n) \otimes (v_1 \otimes \cdots \otimes v_p)) = \\
 \sum_{\sigma \in \mathfrak{S}_{p,n-p}} (-1)^{\epsilon(\sigma)}
 (g(u_{\sigma(1)} \wedge \cdots \wedge u_{\sigma(p)}) \wedge u_{\sigma(p+1)} \wedge \cdots \wedge u_{\sigma(n)}) 
 \otimes (v_1 \otimes \cdots \otimes v_p).
\end{multline}

\bibliography{eduII}

\end{document}